\documentclass[10pt,final,technote]{IEEEtran}

\usepackage[dvips]{graphicx}
\usepackage{float}
\usepackage{times}
\usepackage{amsmath}
\usepackage{amsfonts, amsthm}
\renewcommand{\baselinestretch}{1.14}
\usepackage{cite}

\newtheorem{thm}{Theorem}
\newtheorem{lemma}{Lemma}
\theoremstyle{remark}
\newtheorem*{remark}{Remark}

\theoremstyle{definition}
\newtheorem{exa}{Example}

\DeclareMathOperator{\sech}{sech}
\newcommand{\be}{\begin{equation}}
\newcommand{\ee}{\end{equation}}
\newcommand{\bes}{\begin{equation*}}
\newcommand{\ees}{\end{equation*}}
\newcommand{\bea}{\begin{eqnarray}}
\newcommand{\eea}{\end{eqnarray}}
\newcommand{\bi}{\begin {itemize}}
\newcommand{\ei}{\end {itemize}}
\newcommand{\bmn}{\begin{minipage}}
\newcommand{\emn}{\end{minipage}}
\newcommand{\bfig}{\begin{figure}}
\newcommand{\efig}{\end{figure}}

\newcommand{\bcls}{\begin{columns}}
\newcommand{\ecls}{\end{columns}}
\newcommand{\bcl}{\begin{column}}
\newcommand{\ecl}{\end{column}}

\newcommand{\ptl}{\partial}
\newcommand{\td}{\tilde}

\def\al{\alpha}

\def\v{{\bf v}}
\def\w{{\bf w}}

\def\hv{{\hat{\bf v}}}

\begin{document}

\title{Properties and Applications of a Restricted HR Gradient Operator}

\author{
{Mengdi Jiang{\small $~^{\mathrm{a}}$}, Yi Li{\small $~^{\mathrm{b}}$},
Wei Liu{\small $~^{\mathrm{a}}$}}\vspace{1.6mm}\\
\fontsize{10}{10}\selectfont$^{\mathrm{a}}$\,Communications Research Group, Department of Electronic and Electrical Engineering\\ University of Sheffield, Sheffield, S1 3JD, United Kingdom\\\fontsize{9}{9}\selectfont\ttfamily\upshape\,\{mjiang3, w.liu\}@sheffield.ac.uk\\
\fontsize{10}{10}\selectfont\rmfamily
$^{\mathrm{b}}$\,School of Mathematics and Statistics\\ University of Sheffield, Sheffield, S3 7RH, United Kingdom\\\fontsize{9}{9}\selectfont\ttfamily\upshape\,yili@sheffield.ac.uk\\
}

\maketitle
\begin{abstract}
  For quaternionic signal processing algorithms, the gradients of a
  quaternion-valued function are required for gradient-based methods. Given
  the non-commutativity of quaternion algebra, the definition of the gradients
  is non-trivial. The HR gradient operator provides a viable framework
  and has found a number of applications. However, the applications so far
  have been mainly limited to real-valued quaternion functions and linear
  quaternion-valued   functions. To generalize the operator
  to nonlinear quaternion functions, we define a restricted version of the HR
  operator. The restricted HR gradient operator comes in two versions, the
  left and the right ones. We then present a detailed analysis of the
  properties of the operators, including several different
  product rules and chain rules. Using the new rules, we derive explicit expressions for the derivatives of a
  class of regular nonlinear quaternion-valued functions, and prove that
  the restricted HR gradients are
  consistent with the gradients in real domain.
\end{abstract}
\section{Introduction}
Quaternion calculus has been introduced in signal processing  with  application areas involving three or four-dimensional signals, such as  color image processing~\cite{pei99a,sangwineell2000,parfieniuk10a}, vector-sensor array systems~\cite{BihanN2004,miron06a,lebihan07a,tao13a,tao13b,tao14a,liu14e} and wind profile prediction~\cite{took09a,liu14f}. Several quaternion-valued adaptive
filtering algorithms have been proposed in~\cite{liu14e,liu14f,liu14g,took09a,tao14a}.
Notwithstanding the advantages of the quaternionic algorithms, extra cares
have to be taken in their developments, in particular when the derivatives of
quaternion-valued functions are involved, due to the fact that quaternion
algebra is non-commutative. A so-called HR gradient operator was
proposed in~\cite{mandic2011a}, and has been applied in~\cite{mandic14a}.
The interesting formulation appears to
provide a general and flexible framework that could potentially have wide
applications. However, it has only been applied to real-valued functions
and linear quaternion-valued functions. In order to consider more general
quaternion-valued functions, we propose a pair of restricted
HR gradient operators, the left and the right restricted HR gradient operators,
based on the previous work on the HR gradient operator
\cite{mandic2011a} and our recent work \cite{liu14f}.

To summarize, we make the following main
contributions. Firstly, we
give a detailed derivation of the relation
between the gradients and the increment of a quaternion
function, highlighting the difference between the left and the right gradients
due to the non-commutativity of quaternion algebra.
Secondly, we document several properties of the
operators that have not been reported before, in particular several different versions of
product rules and chain rules.
Thirdly, we derive a general formula for the restricted HR derivatives of a wide class of
regular quaternion-valued nonlinear functions, among which are the exponential, logarithmic, and the
hyperbolic tangent functions. Finally, we prove that the restricted HR
gradients are consistent with the usual definition for the gradient of a real
function of a real variable. Its application to the derivation of quaternion-valued least
mean squares (QLMS) adaptive algorithm is also briefly discussed.

The paper is organised as follows. The restricted HR gradient operator is developed in Section \ref{sec:HRgradient}, with its properties and rules introduced in Section \ref{sec:properties}. Explicit expressions for the derivatives for a wide range
of functions are derived in Section \ref{sec:functions} and results for the right restricted HR operator are summarised in Section \ref{sec:rightHR}. The increment of a general quaternion function is discussed in Section \ref{sec:increment} with the QLMS adaptive algorithm revisited as a special case where the cost function is real-valued. Conclusions are drawn in Section \ref{sec:concl}.

\section{The restricted HR gradient operators}
\label{sec:HRgradient}
\subsection{Introduction of quaternion}

Quaternion is a non-commutative extension of complex number.
A quaternion $q$ is composed  of four parts,
i.e., $q =q_{a} + q_{b}i + q_{c}j + q_{d}k$, where $q_a$ is the real part,
which is also denoted as $R(q)$. The other three terms constitute the
imaginary part $I(q)$.
$i$,
$j$ and $k$ are the three imaginary units, which satisfy the following rules
for multiplication: $ij=k$, $jk=i$, $ki=j$,
$i^{2} = j^{2} = k^{2} = -1$, and
\be \label{eq:ij}
ij=-ji, ki=-ik, kj=-jk.
\ee
Due to (\ref{eq:ij}), in general the product of two quaternions depends on the order,
i.e., $qp \neq pq$ where $p$ and $q$ are quaternions. However, the product
commutes as long as at least one of the factors, say $q$, is real.

Let $v =
|I(q)|$ and $\hv = I(q)/v$, the quaternion $q$ can also be written as $q=q_a +
v \hat{\v}$. $\hv$ is a pure unit quaternion, which has the convenient property
$\hv^2 :=\hv \hv = -1$. The quaternionic conjugate of $q$ is $q^* = q_a - q_b
i - q_c j - q_d k$, or $q^* = q_a - v \hv$. It is easy to show that $qq^* =
q^* q = |q|^2$, and hence $q^{-1} = q^*/|q|^2$.


\subsection{Definition of the restricted HR gradient operators}

Let $f:H\to H$ be a quaternion-valued function of a quaternion $q$,
where $H$ is the non-commutative algebra of quaternions.
We use the notation
$f(q)=f_{a} + f_{b}i + f_{c}j + f_{d}k$,
where $f_a, ..., f_d$ are the components of $f$.
$f$ can also be viewed as
a function of the four components of $q$, i.e., $f = f(q_a, q_b, q_c, q_d)$.
In this view $f$ is a quaternion-valued function on $R^4$: $f: R^4 \to H$.
To express the four real components of $q$, it is
convenient to use its involutions $q^\nu:= -\nu q \nu$ where $\nu \in\{i,
j,k\}$~\cite{ell2007a}. Explicitly, we have
\begin{align}
q^{i}&=-iqi=q_{a} + q_{b}i - q_{c}j - q_{d}k,\\
q^{j}&=-jqj=q_{a} - q_{b}i + q_{c}j - q_{d}k,\\
q^{k}&=-kqk=q_{a} - q_{b}i - q_{c}j + q_{d}k.
\end{align}
The real components can be recovered by
\begin{align}
&q_{a}=\frac{1}{4}(q+q^{i}+q^{j}+q^{k}),
q_{b}=\frac{1}{4i}(q+q^{i}-q^{j}-q^{k}),\label{eq:qab}\\
&q_{c}=\frac{1}{4j}(q-q^{i}+q^{j}-q^{k}),
q_{d}=\frac{1}{4k}(q-q^{i}-q^{j}+q^{k}).\label{eq:qcd}
\end{align}
Two useful relations are
\be
q^* = \frac{1}{2}(q^i + q^j + q^k -  q), \ q + q^i + q^j + q^k = 4 R(q).
\ee
A so-called HR gradient of $f(q)$ was introduced in
~\cite{mandic2011a}, which has been applied to real-valued functions and
linear quaternion-valued functions. In order to find the gradients of
more general quaternion-valued functions, we follow a similar approach to
propose a `restricted' HR gradient operator (some of the derivation was first
presented in \cite{liu14f}).
To motivate the definitions, we consider the differential $df(q)$ with
respect to differential $dq := dq_a  + dq_b i + dq_c j + dq_d k$.
We observe that $df=df_a+idf_b+jdf_c+kdf_d$, where
\be
df_a
=\frac{\partial f_a}{\partial q_a}dq_a+\frac{\partial f_a}{\partial q_{b}}dq_{b}
+\frac{\partial f_a}{\partial q_{c}}dq_{c}+\frac{\partial f_a}{\partial
q_{d}}dq_d.
\ee
We have $dq_a = (dq + dq^i + dq^j + dq^k)/4$ according to (\ref{eq:qab}). Making use
of this and similar expressions for $dq_b$, $dq_c$ and $dq_d$, we
find an expression for $df_a$ in terms of the differentials $dq$, $dq^i$,
$dq^j$ and $dq^k$. Repeating the calculation for $idf_b$, $jdf_c$ and $kdf_d$,
we finally arrive at
\be
df = D dq + D_i dq^i + D_j dq^j + D_k dq^k
\ee
where
\begin{align}
  D&:=\frac{1}{4}\left(\displaystyle\frac{\partial{f}}{\partial q_a}-
\frac{\partial{f}}{\partial q_b} i-
\displaystyle\frac{\partial{f}}{\partial q_c} j-
\displaystyle\frac{\partial{f}}{\partial q_d} k\right), \label{eq:d} \\
  D_i&:=\frac{1}{4}\left(\displaystyle\frac{\partial{f}}{\partial q_a}-
\frac{\partial{f}}{\partial q_b} i+
\displaystyle\frac{\partial{f}}{\partial q_c} j+
\displaystyle\frac{\partial{f}}{\partial q_d} k\right),\label{eq:di}\\
  D_j&:=\frac{1}{4}\left(\displaystyle\frac{\partial{f}}{\partial q_a}+
\frac{\partial{f}}{\partial q_b} i-
\displaystyle\frac{\partial{f}}{\partial q_c} j+
\displaystyle\frac{\partial{f}}{\partial q_d} k\right), \label{eq:dj}\\
  D_k&:=\frac{1}{4}\left(\displaystyle\frac{\partial{f}}{\partial q_a}+
\frac{\partial{f}}{\partial q_b} i+
\displaystyle\frac{\partial{f}}{\partial q_c} j-
\displaystyle\frac{\partial{f}}{\partial q_d} k\right).  \label{eq:dk}
\end{align}
More details are given in Appendix \ref{app:df}.
Thus one may define the partial derivatives of $f(q)$ as follows:
\be\label{eq:dfdq}
\frac{\ptl f}{\ptl q} := D, \quad
\frac{\ptl f}{\ptl q^i} := D_i,
\frac{\ptl f}{\ptl q^j} := D_j,
\frac{\ptl f}{\ptl q^k} := D_k.
\ee
Introducing operators $\nabla_q :=({\partial}/{\partial q},{\partial}/{\partial
q^{i}},{\partial}/{\partial q^{j}},{\partial}/{\partial q^{k}})$,
and $\nabla_r:= ({\partial}/{\partial q_a},{\partial}/{\partial
q_b},{\partial}/{\partial q_c},{\partial}/{\partial q_d})$,
equations (\ref{eq:d}-\ref{eq:dfdq}) may be written as
\begin{equation} \label{eq:hrgrad}
    \nabla_q f =\nabla_r f J^{H}
\end{equation}
where the Jacobian matrix
\begin{equation}
    J=\frac{1}{4}
    \begin{bmatrix}
        1 & i & j & k
        \\
        1 & i & -j & -k
        \\
        1  & -i & j & -k
         \\
        1  & -i & -j & k
      \end{bmatrix}
\end{equation}
and $J^H$ is the Hermitian transpose of $J$ \cite{mandic2011a}.
Using $JJ^H = J^H J= 1/4$   \cite{mandic14a},
we may also write
\begin{equation}
    \nabla_q f J=\frac{1}{4}\nabla_r f,
\end{equation}
which is the inverse formulae for the derivatives.

We call the gradient operator defined by (\ref{eq:hrgrad}) the restricted HR
gradient operator. The operator is closely related to the HR operator
introduced in~\cite{mandic2011a}. However, in the original definition of the
HR operator, the Jacobian $J$ appears on the left-hand side of $\nabla_r f$,
whereas in our definition it appears on the right (as the Hermitian transpose).

The differential $df$ is related to $\nabla_q f$ by
\be \label{eq:df}
df = \frac{\ptl f}{\ptl q} dq
+ \frac{\ptl f}{\ptl q^i} dq^i
+ \frac{\ptl f}{\ptl q^j} dq^j
+\frac{\ptl f}{\ptl q^k} dq^k .
\ee
Due to the non-commutativity of quaternion products, the order of the
 factors in the products of the above equation (as well as
equations (\ref{eq:d}-\ref{eq:dk})) can not be swapped. In fact, one may call the
    above operator the left restricted HR gradient operator. As is shown in Appendix
    \ref{app:df}, one may also define a right restricted HR gradient operator by
    \be \label{eq:rgradop}
    (\nabla^R_q f)^T := J^* (\nabla_r f)^T,
    \ee
where
$$\nabla^R_q:=(\ptl^R/\ptl q, \ptl^R/\ptl q^i, \ptl^R/\ptl
q^j, \ptl^R/\ptl q^k),$$
and
\begin{align}
  \frac{\ptl^R f}{\ptl
  q}&:=\frac{1}{4}\left(\displaystyle\frac{\partial{f}}{\partial q_a}- i
\frac{\partial{f}}{\partial q_b} -
j\displaystyle\frac{\partial{f}}{\partial q_c} -
k\displaystyle\frac{\partial{f}}{\partial q_d} \right), \label{eq:rdfdq} \\
\frac{\ptl^R f}{\ptl q^i}&:=\frac{1}{4}\left(\displaystyle\frac{\partial{f}}{\partial q_a}-
i\frac{\partial{f}}{\partial q_b} +
j\displaystyle\frac{\partial{f}}{\partial q_c} +
k\displaystyle\frac{\partial{f}}{\partial q_d} \right),\label{eq:rdfdqi}\\
\frac{\ptl^R f}{\ptl q^j}&:=\frac{1}{4}\left(\displaystyle\frac{\partial{f}}{\partial q_a}+
i\frac{\partial{f}}{\partial q_b} -
j\displaystyle\frac{\partial{f}}{\partial q_c} +
k\displaystyle\frac{\partial{f}}{\partial q_d} \right), \label{eq:rdfdqj}\\
\frac{\ptl^R f}{\ptl q^k}  &:=\frac{1}{4}\left(\displaystyle\frac{\partial{f}}{\partial q_a}+
i\frac{\partial{f}}{\partial q_b} +
j\displaystyle\frac{\partial{f}}{\partial q_c} -
k\displaystyle\frac{\partial{f}}{\partial q_d} \right).  \label{eq:rdfdqk}
\end{align}
The right restricted HR gradient operator is
    related to the differential $df$ by
    \be
    df = dq\frac{\ptl^R f}{\ptl q} + dq^i  \frac{\ptl^R f}{\ptl q^i} + dq^j
    \frac{\ptl^R f}{\ptl q^j} + dq^k \frac{\ptl^R f}{\ptl q^k}.
    \ee
    In general, the left and right restricted HR gradients are not the same.
    For example, even for the simplest linear function
    $f(q) = q_0 q$ with $q_0 \in H$ a constant, we have
    \be
    \frac{\ptl q_0 q}{\ptl q} = q_0, \ \frac{\ptl^R q_0 q} {\ptl q} = R(q_0).
    \ee
    However, we will show later that the two gradients coincide for a class of functions.
    In particular, they are the same for real-valued quaternion functions.

The relation between the gradients and the differential is an important ingredient of
gradient-based methods, which we will discuss further later.

\section{Properties and rules of the operator}
\label{sec:properties}

We will now focus on the left restricted HR gradient and simply call it the
restricted HR gradient unless stated otherwise.
It can be easily calculated from the definitions, that
\begin{equation}
  \frac{\partial q}{\partial q}=1,~\frac{\ptl q^\nu}{\ptl q}
  =0,~\frac{\partial q^*}{\partial q}=-\frac{1}{2},
 \end{equation}
 where $\nu \in \{i,j,k\}$.
However, in order to find the derivatives for more complex quaternion
functions, it is useful to first establish the rules of the gradient
operators. We will see that some of the usual
rules do not apply due to the non-commutativity of quaternion products.
\begin{enumerate}
  \item Left-linearity: for arbitrary constant quaternions $\alpha$ and
    $\beta$, and functions $f(q)$ and $g(q)$, we have
    \begin{equation}
      \frac{\partial (\alpha f + \beta g)}{\partial q^\nu} = \alpha \frac{\partial
      f}{\partial q^\nu} + \beta \frac{\partial g}{\partial q^\nu}
    \end{equation}
  for $\nu \in \{1,i,j,k\}$ with $q^1:=q$.
  However, linearity does not hold for right multiplications, i.e., in general
    \begin{equation}
      \frac{\partial f\alpha} {\partial q} \neq \frac{\partial f}{\partial q}
      \alpha.
    \end{equation}
    This is because, according to the definition (\ref{eq:d}),
    \be
    \frac{\ptl f \al}{\ptl q} = \frac{1}{4}\sum_{(\phi,\gamma)}\frac{\ptl
    f}{\ptl q_\phi} \al \gamma
    \ee
    for $(\phi, \gamma) \in \{(a,1), (b, -i), (c, -j), (d, -k)\}$.
    However, $\al \gamma \neq \gamma\al$ in general. Therefore it is different from
    $(\ptl f/\ptl q) \al$, which is
    \be
    \frac{1}{4}\left(\frac{\ptl f}{\ptl q_a}  -
\frac{\ptl f}{\ptl q_b} i
-\frac{\ptl f}{\ptl q_c}  j
-\frac{\ptl f}{\ptl q_d}  k
    \right)\al .
    \ee

 \item The first product rule: the following product rule holds:
   \begin{equation} \label{eq:prodq}
   \nabla_q(fg) = f \nabla_q g+[(\nabla_r f)g]J^{H}.
  \end{equation}
  For example,
\be
\frac{\ptl fq}{\ptl q} = f \frac{\ptl g}{\ptl q} + \frac{1}{4}\left(\frac{\ptl
  f}{\ptl q_a}g - \frac{\ptl f}{\ptl q_b}gi - \frac{\ptl f}{\ptl q_c}gj -
\frac{\ptl f}{\ptl q_d}gk\right).
\ee
  Thus the product rule in general is different from the usual one.
  \item The second product rule: However, the usual product rule applies to differentiation with respect to real variables, i.e.,
    \begin{equation}
      \frac{\partial fg}{\partial q_{\phi}} =\frac{\partial f}{\partial q_{\phi}} g + f\frac{\partial g}{\partial q_{\phi}}
    \end{equation}
    for $\phi = a, b, c,$ or $d$.
  \item The third product rule: The usual product rule also applies if at least one of the two
    functions $f(q)$ and $g(q)$ is real-valued, i.e.,
    \be
    \frac{\ptl fq}{\ptl q} = f \frac{\ptl g}{\ptl q} + \frac{\ptl f}{\ptl q}
    g.
    \ee
  \item The first chain rule:
   For a composite function $f(g(q))$, $g(q):=g_a + g_b i + g_c j + g_d k$ being a quaternion-valued
   function, we have the following chain rule \cite{mandic14a}:
\begin{equation}
\nabla_q f = (\nabla ^{g}_{q}f) M
\end{equation}
where
  $\nabla ^{g}_{q}:=({\partial }/{\partial g},{\partial }/{\partial g^i},
  {\partial }/{\partial g^j},{\partial }/{\partial g^k})$ and $M$ is a
  $4\times 4$ matrix with element $M_{\mu \nu} =
  \ptl g^\mu/\ptl q^\nu$ for $\mu, \nu \in \{1,i,j,k\}$ and $g^\mu = - \mu g \mu$
  ($g^1$ is understood as the same as $g$). Explicitly, we may write
  \be
  \frac{\ptl f}{\ptl q^\nu} = \sum_{\mu} \frac{\ptl f}{\ptl g^\mu} \frac{\ptl
  g^\mu}{\ptl q^\nu}.
  \ee
  The proof is outlined in Appendix \ref{app:chain}.
\item The second chain rule:
  The above chain rule uses $g$ and its involutions as the intermediate
  variables. It is sometimes convenient to use the real components of $g$
  for that purpose instead. In this case, the following chain rule may be used:
\be
   \nabla_q f  =(\nabla ^{g}_{r}f)O
\ee
where $O$ is a $4\times 4$ matrix with entry $O_{\phi \nu} = \ptl g_\phi/\ptl
q^\nu$ with $\phi\in \{a, b, c, d\}$ and $\nu \in \{1,i,j,k\}$, and
$\nabla_r^g :=(\ptl/\ptl g_a, \ptl/\ptl g_b, \ptl/\ptl g_c, \ptl /\ptl g_d)$.
Explicitly, we have
\begin{equation}
\frac{\partial f}{\partial q^{\nu}}
=\sum_{\phi}\frac{\partial f}{\partial g_\phi}\frac{\partial g_\phi}{\partial
q^{\nu}}.
\end{equation}

\item The third chain rule: if the intermediate function $g(q)$ is
  real-valued, i.e., $g=g_a$, then from the second chain rule, we obtain
  \be \label{eq:chain3}
  \frac{\ptl f}{\ptl q^\nu} = \frac{\ptl f}{\ptl g} \frac{\ptl g}{\ptl q^\nu}.
  \ee

  \item $f(q)$ is not independent of $q^i$, $q^j$ or $q^k$ in the sense
      that, in general,
    \begin{equation}
      \frac{\partial f(q)}{\partial q^i}  \neq 0,
      \frac{\partial f(q)}{\partial q^j}  \neq 0,
      \frac{\partial f(q)}{\partial q^k}  \neq 0.
    \end{equation}
    This can be illustrated by $f(q) = q^2$. Using the first product rule (equation
  (\ref{eq:prodq})), we have
    $$
    \frac{\ptl q^2}{\ptl q^i} = q\frac{\ptl q}{\ptl q^i} + \frac{1}{4}
    \sum_{(\phi, \nu)}
    \frac{\ptl q}{\ptl q_\phi} q \nu
    $$
    for $(\phi, \nu) \in \{(a,1), (b, i), (c, -j), (d, -k)\}$.
    It can then be shown that
    \be \label{eq:q2}
    \frac{\partial q^2}{\partial q^i} = q_b i, ~
    \dfrac{\partial q^{2}}{\partial q^{j}}=q_c j,~ \dfrac{\partial
    q^{2}}{\partial q^{k}}=q_d k.
    \ee
    This property demonstrates the intriguing difference between
    the HR derivative and the usual derivatives, although we can indeed show that
    \be
    \frac{\ptl q}{\ptl q^{\nu}} = 0.
    \ee
    One implication of this observation is that, for a nonlinear algorithm involving
    simultaneously more than one gradients $\ptl f/\ptl q^\nu$, we have to
    take care to include all the terms.

\end{enumerate}

\section{Restricted HR derivatives for a class of regular functions}
\label{sec:functions}
Using the above operation rules, we may find  explicit expressions for the derivatives for a whole range
of functions. We first introduce the following lemma:

\begin{lemma} \label{lm:power}
  The derivative of the power function $f(q) = (q-q_0)^n$, with integer $n$ and
  constant quaternion $q_0$,
  is
    \begin{equation} \label{eq:power}
      \frac{\partial f(q)}{\partial q} = \frac{1}{2} \left(n
        \td{q}^{n-1} + \frac{\td{q}^{n} -
      \td{q}^{*n}}{\td{q} - \td{q}^*}\right) ,
    \end{equation}
    with $\td{q} = q-q_0$.
\end{lemma}
\begin{remark}
  The division in $(\td{q}^{n} - \td{q}^{*n})/(\td{q}-\td{q}^*)$ is understood as
  $(\td{q}^n-\td{q}^{*n})(\td{q}-\td{q}^*)^{-1}$ or
  $(\td{q}-\td{q}^*)^{-1}(\td{q}^n-\td{q}^{*n})$ which are the
    same since the two factors commute.  The division operations in what
    follows are understood in the same way.
  \end{remark}

    \begin{IEEEproof}
    The lemma is obviously true for $n=0$.
    Let $n\ge 1$, we
    apply the first product rule, and find
    \begin{equation}
      \frac{\partial (q-q_0)^n}{\partial q}=\td{q}\frac{\partial
      \td{q}^{n-1}}{\partial q}+R(\td{q}^{n-1})
    \end{equation}
    where $R(\td{q}^{n-1})$ is the real part of $\td{q}^{n-1}$. We then obtain
    by induction
    \begin{equation}
      \frac{\partial (q-q_0)^n}{\partial q}=\sum^{n-1}_{m=0}\td{q}^{m}R(\td{q}^{n-1-m}).
    \end{equation}
    Using $R(\td{q}^{n-1-m})=\frac{1}{2}(\td{q}^{n-1-m}+\td{q}^{*(n-1-m)})$, the summations
    can be evaluated explicitly, leading to
    equation (\ref{eq:power}).

    For $n<0$, we use the recurrent relation
    \begin{equation}
      \frac{\partial ((q-q_0)^{-n})}{\partial q}=\td{q}^{-1}\left[\frac{\ptl
      \td{q}^{-(n-1)}}{\ptl q} -R(\td{q}^{-n})\right]
    \end{equation}
    and the result
    \be
    \frac{\ptl (q-q_0)^{-1}}{\ptl q} = - \td{q}^{-1}R(\td{q}^{-1}).
    \ee
    Equation (\ref{eq:power}) is proven by using induction as for $n>0$.  More
    details are given in Appendix \ref{app:lm1}.
    \end{IEEEproof}

    \begin{thm} \label{thm:dfdq}
      Assuming $f: H\to H$ admits a power series representation $f(q) :=
      g(\td{q}) :=
  \sum_{n=-\infty}^\infty a_n \td{q}^n$, with $a_n$ being a quaternion constant and
  $\td{q} = q- q_0$, for $R_1\le |\td{q}|\le R_2$ with $R_1, R_2>0$ being some
  constants, then
  \be
  \frac{\ptl f(q)}{\ptl q} = \frac{1}{2}\left[f'(q) + (g(\td{q}) -
  g(\td{q}^*))(\td{q}-\td{q}^*)^{-1}\right],
  \ee
  where $f'(q)$ is the derivative in the usual sense, i.e.,
  \be
f'(q):=\sum_{n=-\infty}^\infty n a_n \td{q}^{n-1}=\sum_{n=-\infty}^\infty n a_n (q-q_0)^{n-1}.
  \ee
\end{thm}
\begin{IEEEproof} Using Lemma \ref{lm:power} and the left-linearity of HR
  gradients, we have
  \begin{align}
    \frac{\ptl f}{\ptl q}& = \frac{1}{2}\sum_{n=-\infty}^\infty a_n [n \td{q}^{n-1}
    + (\td{q}^n
  -\td{q}^{*n})(\td{q}-\td{q}^*)^{-1}] \notag\\
  & = f'(q)+ \frac{1}{2}
  \left[\sum_{n=\infty}^\infty a_n(\td{q}^n - \td{q}^{*n})\right](\td{q}-\td{q}^*)^{-1}\notag\\
  & = \frac{1}{2} [f'(q) + (g(\td{q})-g(\td{q}^*))(\td{q}-\td{q}^*)^{-1}],\notag
\end{align}
proving the theorem.
\end{IEEEproof}
The functions $f(q)$ form a class of regular functions on $H$. A
full discussion of such functions is beyond the scope of this
paper.
However, we
  note that a similar class of functions have been discussed in
  \cite{GentiliStruppa07}. A parallel development for the former is possible,
  and will be the topic of a future paper. Meanwhile, we observe that
  many useful elementary functions satisfy the conditions in Theorem
  \ref{thm:dfdq}. To
illustrate the application of the theorem, we list below the derivatives of a
number of such functions.
\begin{exa}
  Exponential function $f(q) = e^q$ has representation
  \be
  e^q:=\sum_{n=0}^\infty \frac{q^n}{n!}.
  \ee
  Applying Theorem \ref{thm:dfdq} with $a_n = 1/n!$ and $q_0 = 0$, we have
    \begin{equation} \label{eq:exp}
      \frac{\partial e^q}{\partial q} = \frac{1}{2} \left(e^q + \frac{e^{q} -
      e^{q*}}{q - q^*}\right).
    \end{equation}
     Making use of $e^q = e^{q_a} (\cos v + \hv \sin v)$ and $q= q_a+ \hv v$,
     we have
    \begin{equation}
      \frac{\partial e^q}{\partial q} = \frac{1}{2} \left(e^q +e^{q_a}
      v^{-1}\sin v \right).
    \end{equation}
\end{exa}
\begin{exa}
  The logarithmic function $f(q) = \ln q $ has representation
  \be
  \ln q = \sum_{n=1}^\infty \frac{(-1)^{n-1}}{n}(q-1)^n.
  \ee
  with $a_n = (-1)^{n-1}/n$ and $q_0 =1$. Since $q_0$ is a real number,
  $g(\td{q}^*) = f(q^*)$.
  Therefore, we have from Theorem \ref{thm:dfdq}
    \be \label{eq:ln1}
    \frac{\ptl \ln q}{\ptl
    q}=\frac{1}{2}\left(q^{-1}+\frac{\ln q - \ln q^*}{q-q^*}\right).
    \ee
    Using representation $ \ln q=\ln {|q|}+\hv \arccos ({q_a}/{|q|})$, the
    expression can be simplified as
    \be \label{eq:ln}
    \frac{\ptl \ln q}{\ptl
    q}=\frac{1}{2}\left(q^{-1}+\frac{1}{v}\arccos{\frac{q_a}{|q|}}\right),
    \ee
    where $v = |I(q)|$.
\end{exa}
\begin{exa}
Hyperbolic tangent function $f(q) = \tanh q$ is defined as
\be
\tanh q :=\frac{e^q - e^{-q}}{e^q + e^{-q}} = q - \frac{q^3}{3} + \frac{2
q^5}{15} - ...
\ee
Therefore, Theorem \ref{thm:dfdq} applies.
     On the other hand, using the relation $e^q = e^{q_a} (\cos v + \hv \sin v)$, we can
     show that
     \be
     \tanh q =
     \frac{1}{2}\frac{\sinh2q_a+\hv\sin2v}{\sinh^2 q_a+\cos^2 {v}}.
     \ee
     Then the second term in the expression given by Theorem
     \ref{thm:dfdq} can be simplified. The final expression can be written as
     \begin{equation}
       \dfrac{\partial \tanh q}{\partial q}=\frac{1}{2}\left(\sech^2
       q+\frac{v^{-1}\sin2v}{\cosh2q_a+\cos2v}\right),
     \end{equation}
     where $\sech q := 1/\cosh q$ is the quaternionic hyperbolic secant function.
\end{exa}
\begin{remark}
Apparently,
the derivatives for these functions can also be
  found by direct calculations without resorting to Theorem \ref{thm:dfdq}.
\end{remark}

  We now turn to a question of more theoretical interests.
  Even though it might not be obvious from the definitions,
    the following theorem shows that the restricted HR derivative is
    consistent with the derivative in the real domain for a class of
    functions, including those in the above examples.
\begin{thm} \label{thm:consistency} For the function $f(q)$ in Theorem
  \ref{thm:dfdq}, if $q_0$ is a real number, then
  \be
  \frac{\ptl f(q)}{\ptl q} \to f'(q)
  \ee
  when $q \to R(q) $, i.e., when $q$ approaches a real number.
\end{thm}
\begin{IEEEproof} Using the polar representation, we write $\td{q} = |\td{q}| \exp(\hv \theta)$, where
  $\theta = \arcsin (v/|\td{q}|)$ is the argument of $\td{q}$ with $v =
  |I(\td{q})|$.  Then $\td{q}^n = |\td{q}|^n \exp(n \hv \theta)$, and
  \be
  (\td{q}^n - \td{q}^{*n})(\td{q}-\td{q}^*)^{-1} = \frac{
  I(\td{q}^n)}{I(\td{q})} =
  \frac{|\td{q}|^{n-1} \sin (n \theta)}{\sin\theta}.
  \ee
  For real $q_0$, $\td{q}\to q_a - q_0$ and $v\to 0$ when $q\to
  R(q)$. As a consequence,
  $\theta \to 0$ at the limit (or $\theta \to \pi$, which can be dealt with by
  slight modification), and
\be
\frac{\sin (n \theta)}{\sin\theta} \sim \frac{\sin (n \theta)}{\theta} \to n,
\quad
|\td{q}|^{n-1} \to (q_a - q_0)^{n-1}.
\ee
Therefore,
\be
(\td{q}^n - \td{q}^{*n})(\td{q}-\td{q}^*)^{-1} \to n \td{q}^{n-1}
\ee
and
\be
[g(\td{q})-g(\td{q}^*)](\td{q}-\td{q}^*)^{-1} \to
\sum_{n=-\infty}^{\infty} n a_n \td{q}^{n-1} = f'(q).
\ee
Thus
\be
\frac{\ptl f(q)}{\ptl q} \to \frac{1}{2}[f'(q) + f'(q)] = f'(q).
\ee
\end{IEEEproof}
The functions in above three examples all satisfy the conditions in Theorem
\ref{thm:consistency}, hence we expect Theorem \ref{thm:consistency} applies.
One can easily verify by direct calculations that the theorem indeed holds.

\section{The right restricted HR gradients}
\label{sec:rightHR}
In this section, we briefly summarize the results for the right restricted HR
gradients, and highlight the difference with left restricted HR gradients.
\begin{enumerate}
  \item Right-linearity: for arbitrary quaternion constants $\alpha$ and
    $\beta$, and functions $f(q)$ and $g(q)$, we have
    \begin{equation}
      \frac{\partial^R (f \alpha + g\beta )}{\partial q^\nu} = \frac{\partial^R
      f}{\partial q^\nu}\alpha  + \frac{\partial^R g}{\partial q^\nu}\beta.
    \end{equation}
  However, linearity does not hold for left multiplications, i.e., in general
    \begin{equation}
      \frac{\partial^R \alpha f} {\partial q} \neq \alpha\frac{\partial^R f}{\partial q}.
    \end{equation}
 \item The first product rule: for the right restricted HR operator, the following product rule holds:
   \begin{equation} \label{eq:rprodq}
     [\nabla^R_q(fg)]^T =  [(\nabla^R_q f) g]^T+J^{*}[f(\nabla_r g)^T].
  \end{equation}
  The second and third product rules are the same as for the left restricted
  operator.
\item The first chain rule: for the composite function $f(g(q))$, we have
  \be
(\nabla^R_q f)^T = M^T(\nabla ^{gR}_{q}f)^T.
\ee
\item The second chain rule becomes:
  \be
  (\nabla_q^R f)^T = O^T (\nabla_r^g f)^T.
  \ee
\item The third chain rule becomes
  \be \label{eq:rchain3}
  \frac{\ptl^R f}{\ptl q^\nu} = \frac{\ptl g}{\ptl q^\nu}\frac{\ptl f}{\ptl g} .
  \ee
  Note that, $\ptl g/\ptl q^\nu = \ptl^R g/\ptl q^\nu$ since $g$ is
  real-valued.
  We thus have omitted the superscript $R$. Also, $\ptl f/\ptl g$ is a
  real derivative, so there is no distinction between left and right
  derivatives.
\end{enumerate}
We can also find the right restricted HR gradients for common quaternion
functions. First of all, Lemma \ref{lm:power} is also true for right
derivatives:
\begin{lemma}\label{lm:rpower} For $f(q) = (q-q_0)^n$ with $n$ integer and $q_0$ a constant
  quaternion, we have
    \begin{equation} \label{eq:rpower}
      \frac{\partial^R f(q)}{\partial q} = \frac{1}{2} \left(n
        \td{q}^{n-1} + \frac{\td{q}^{n} -
      \td{q}^{*n}}{\td{q} - \td{q}^*}\right) ,
    \end{equation}
    with $\td{q} = q-q_0$.
\end{lemma}
\begin{remark} To prove the lemma, we use the following recurrent relations:
    \begin{equation}
      \frac{\partial (q-q_0)^n}{\partial q}=\frac{\partial
      \td{q}^{n-1}}{\partial q}\td{q}+R(\td{q}^{n-1})
    \end{equation}
    \begin{equation}
      \frac{\partial ((q-q_0)^{-n})}{\partial q}=\left[\frac{\ptl
      \td{q}^{-(n-1)}}{\ptl q} -R(\td{q}^{-n})\right]\td{q}^{-1}.
    \end{equation}
\end{remark}
Using Lemma \ref{lm:rpower},
We can prove the following result:
\begin{thm} \label{thm:rdfdq}
      Assuming $f: H\to H$ admits a power series representation $f(q) :=
      g(\td{q}) :=
  \sum_{n=-\infty}^\infty \td{q}^na_n $, with $a_n$ being a quaternion constant and
  $\td{q} = q- q_0$, for $R_1\le |\td{q}|\le R_2$ with $R_1, R_2>0$ being some
  constants, then
  \be
  \frac{\ptl^R f(q)}{\ptl q} = \frac{1}{2}\left[f'(q) + (\td{q}-\td{q}^*)^{-1}(g(\td{q}) -
  g(\td{q}^*))\right],
  \ee
  where $f'(q)$ is the derivative in the usual sense, i.e.,
  \be
f'(q):=\sum_{n=-\infty}^\infty n \td{q}^{n-1}a_n =\sum_{n=-\infty}^\infty n (q-q_0)^{n-1}a_n .
  \ee
\end{thm}
Note that, the functions $f(q)$ in Theorem \ref{thm:rdfdq} in general
form a different class of functions than the one in Theorem
\ref{thm:dfdq}, because in the series representation $a_n$ appears on the
right-hand side of the powers. However, if $a_n$ is a real number, then the
two classes of functions coincide. Therefore, we have the following result:
\begin{thm}
  If $a_n$ is real, then the left and right restricted HR gradients of $f(q)$
  coincide.
\end{thm}
\begin{remark} As a consequence, we can see immediately the right derivatives
  for the exponential, logarithmic and hyperbolic tangent functions are the
  same as the left ones.
\end{remark}
Apparently, Theorem \ref{thm:consistency} is also true for the right
derivatives. Hence, we have:
\begin{thm}
 The right-restricted HR gradient is consistent with the real gradient in the
 sense of Theorem \ref{thm:consistency}.
\end{thm}

\section{The increment of a quaternion function \label{sec:increment}}
When $f(q)$ is a real-valued quaternion function, both left and right restricted HR
gradients are coincident with the HR gradients. Besides,
we have
    \be \label{eq:realdfdq}
    \frac{\ptl^R f}{\ptl q^\nu} = \frac{\ptl f}{\ptl q^\nu} = \left(\frac{\ptl f}{\ptl q}\right)^\nu,
    \ee
    where $\nu \in {i,j,k}$.
    Thus only $\ptl f/\ptl q$ is independent. As a consequence (see also \cite{mandic2011a}),
\begin{align}
df & = \sum_{\nu} \frac{\ptl f}{\ptl q^\nu} d q^\nu =
\sum_{\nu} \left(\frac{\ptl f}{\ptl q}\right)^\nu d q^\nu \notag \\
&=
\sum_{\nu} \left(\frac{\ptl f}{\ptl q} d q \right)^\nu= 4 R\left(\frac{\ptl
f}{\ptl q} dq\right),
\end{align}
where equation (\ref{eq:realdfdq}) has been used.
Hence, $-(\ptl f/\ptl q)^*$
gives the steepest descent direction for $f$, and the increment is determined
by $\ptl f/\ptl q$.

On the other hand, if $f$ is a quaternion-valued function, the increment will
depend on all four derivatives. Taking $f(q) = q^2$ as an example, we have (see
equations (\ref{eq:q2}) and (\ref{eq:power}))
\be
dq^2 = (q + q_a) dq + q_b i dq^i + q_c j dq^j + q_d k dq^k,
\ee
even though $f(q)$ appears to be independent of $q^i$, $q^j$ and
$q^k$. It can be verified that the above expression is the same as the
differential form given in terms of $dq_a$, $dq_b$, $dq_c$ and $dq_d$. Thus it
is essential to include the contributions from $\ptl f/\ptl q^i$ etc.

We also note that, if the right gradient is used consistently, the same
increment would result, since the basis of the definitions is the same,
namely, the differential form in term of $dq_a$, $dq_b$, $dq_c$ and $dq_d$.

Now we apply the quaternion-valued restricted HR gradient operator to develop the QLMS algorithm as an
application. This version of QLMS has been derived in \cite{tao14a,liu14f,liu14g,mandic14a}. However,
with the rules we have derived, some of the calculations can be simplified, as
we will be showing below.

In terms of a standard adaptive filter, the output $y[n]$ and error $e[n]$ can be expressed as
\begin{align}
y[n]&={\textbf{w}^{T}[n]}{\textbf{x}[n]}\\
e[n]&=d[n]-{\textbf{w}^{T}[n]}{\textbf{x}[n]},
\end{align}
where $\textbf{w}[n]=[w[1], w[2], \cdots, w[M]]^{T}$ is the quaternion adaptive weight coefficient
vector with length $M$,
$d[n]$ the reference signal, and
$\textbf{x}[n]=[x[n-1], x[n-2], \cdots, x[n-M]]^{T}$ the quaternion input sample
sequence. The conjugate $\textbf{e}^{*}[n]$ of the error signal $e[n]$ is
\begin{equation}
e^{*}[n]=d^{*}[n]-{\textbf{x}^{H}[n]}{\textbf{w}^{*}[n]}.
\end{equation}
The cost function
is defined as $J[n]=e[n]e^{*}[n]$ which is real-valued.
According to the discussion above and  \cite{mandic2011a,brandwood83a},
the conjugate gradient $(\nabla_{\textbf{w}}J[n])^*$
gives the maximum steepness direction for the optimization surface. Therefore
it is used to update the
weight vector. Specifically,
\be
  \textbf{w}[n+1] = \textbf{w}[n]-\mu (\nabla_{\textbf{w}}J[n])^{*},
\ee
where $\mu$ is the step size. To find $\nabla_\w J$ , we use the first product rule:
\begin{align}
\nabla_\w &=\frac{\partial e[n]e^*[n]}{\partial \textbf{w}}\nonumber\\
&=e[n]\frac{\partial e^*[n]}{\partial \textbf{w}}+\frac{1}{4}(\frac{\partial e[n]}{\partial \textbf{w}_a}e^*[n]-\frac{\partial e[n]}{\partial \textbf{w}_b}e^*[n]i\nonumber\\
&-\frac{\partial e[n]}{\partial \textbf{w}_c}e^*[n]j-\frac{\partial e[n]}{\partial \textbf{w}_d}e^*[n]k)
\end{align}

After some algebra,
we find
\be
\nabla_{\textbf{w}}J[n]=-\frac{1}{2}\textbf{x}[n]e^*[n].
\ee
Therefore we obtain the following update equation for the QLMS algorithm with a step size $\mu$
\be
\textbf{w}[n+1] = \textbf{w}[n]+\mu(e[n]\textbf{x}^{*}[n]).
\label{eq:update_weight_vector}
\ee
Some simulation results have been reported in \cite{liu14f}.

\section{Conclusions}
\label{sec:concl}
We have proposed a restricted HR gradient operator and
discussed its properties, in particular several different versions of
product rules and chain rules. Using the operator, we apply the rules to find
the derivatives for a wide class of nonlinear quaternion-valued functions that
admit a power series representation. The class includes the common elementary
functions such as the exponential function, the logarithmic function, among
others. The explicit expressions for the
derivatives will be useful for nonlinear signal processing applications. We also prove for a
wide class of functions, that the restricted HR gradient tends to the derivatives for real functions with
respect to real variables, when the independent quaternion variable tends to
the real axis, thus showing the consistency of the definition.


\appendices
\section{Definition of the operators \label{app:df}}

We consider $df = df_a + i df_b + j df_c + k df_d$. By definition, we have
$df_\gamma = \sum_{\phi } ({\ptl f_\gamma}/{\ptl q_\phi}) d q_\phi$, with
$\gamma, \phi \in \{a,b,c,d\}$. Using the relations
\begin{align}
  dq_{a}&=\frac{1}{4}(dq+dq^{i}+dq^{j}+dq^{k}),\label{eq:qa}\\
  dq_b&=\frac{1}{4i}(dq+dq^{i}-dq^{j}-dq^{k}),\label{eq:qb}\\
  dq_{c}&=\frac{1}{4j}(dq-dq^{i}+dq^{j}-dq^{k}),\label{eq:qc}\\
  dq_{d}&=\frac{1}{4k}(dq-dq^{i}-dq^{j}+dq^{k}),\label{eq:qd}
\end{align}
we may rewrite $df_\gamma$ as follows
\begin{align}
df_\gamma &=\phantom{+}\frac{1}{4}(\frac{\partial f_\gamma}{\partial q_a}-
  i\frac{\partial f_\gamma}{\partial q_{b}}-j\frac{\partial f_\gamma}{\partial q_{c}}-k\frac{\partial f_\gamma}{\partial q_{d}})dq\nonumber\\
&\phantom{=}+\frac{1}{4}(\frac{\partial f_\gamma}{\partial q_a}-i\frac{\partial f_\gamma}{\partial q_{b}}+j\frac{\partial f_\gamma}{\partial q_{c}}+k\frac{\partial f_\gamma}{\partial q_{d}})dq^{i}\nonumber\\
&\phantom{=}+\frac{1}{4}(\frac{\partial f_\gamma}{\partial q_a}+i\frac{\partial f_\gamma}{\partial q_{b}}-j\frac{\partial f_\gamma}{\partial q_{c}}+k\frac{\partial f_\gamma}{\partial q_{d}})dq^{j}\nonumber\\
&\phantom{=}+\frac{1}{4}(\frac{\partial f_\gamma}{\partial q_a}+i\frac{\partial
f_\gamma}{\partial q_{b}}+j\frac{\partial f_\gamma}{\partial q_{c}}-k\frac{\partial
f_\gamma}{\partial q_{d}})dq^{k}\notag
\end{align}
which can be written as
\be \label{eq:dfgammaapp}
df_\gamma =\frac{1}{4}\sum_{\nu} \left( \sum_{(\phi,\mu)} \frac{\ptl f_\gamma}{\ptl q_\phi}
  \mu^\nu\right) dq^\nu
\ee
where $(\phi,\mu) \in \{(a,1),(b,-i),(c,-j),(d,-k)\}$, $\nu \in \{1,i,j,k\}$,
and $\mu^\nu$ is the $\nu$-involution of $\mu$. Therefore
\begin{align}
  df &= df_a + i df_b + j df_c + k df_d \notag \\
     & = \frac{1}{4}\sum_{\nu} \left( \sum_{(\phi,\mu)} \frac{\ptl (f_a + if_b
+ j f_c + k f_d)}{\ptl q_\phi}
  \mu^\nu\right) dq^\nu \notag \\
  & = \frac{1}{4}\sum_{\nu} \left( \sum_{(\phi,\mu)} \frac{\ptl f}{\ptl q_\phi}
\mu^\nu\right) dq^\nu \label{eq:dfapp}
\end{align}
which leads to the definitions (\ref{eq:d}-\ref{eq:df}) in the main text.
Note that, because $\mu^\nu$ and $dq^\nu$ are quaternions, to obtain the last
equation, we need to multiply $df_b$, $df_c$ and $df_d$ by $i$, $j$, and $k$ from
the left.

On the other hand, we notice that the prefactors in
(\ref{eq:qb}-\ref{eq:qd}) may be moved to the right-hand side of the other
factors, i.e.,
we may write
\begin{align}
  dq_{a}&=(dq+dq^{i}+dq^{j}+dq^{k})\frac{1}{4},\label{eq:qar}\\
  dq_b&=(dq+dq^{i}-dq^{j}-dq^{k})\frac{1}{4i},\label{eq:qbr}\\
  dq_{c}&=(dq-dq^{i}+dq^{j}-dq^{k})\frac{1}{4j},\label{eq:qcr}\\
  dq_{d}&=(dq-dq^{i}-dq^{j}+dq^{k})\frac{1}{4k}.\label{eq:qdr}
\end{align}
Using these relations, we may find another expression for $df_\gamma$
following the procedure above:\be
df_\gamma =\frac{1}{4}\sum_{\nu} dq^\nu \left( \sum_{(\phi,\mu)}  \mu^\nu\frac{\ptl f_\gamma}{\ptl q_\phi}
 \right).
\ee
The expression is
different from (\ref{eq:dfgammaapp}),
in that the differentials $dq^\nu$ are on the left of $\mu^\nu$. Therefore, we
derive
\begin{align}
  df &= df_a +  df_b i+  df_c j+  df_d k\notag \\
     & = \frac{1}{4}\sum_{\nu}  dq^\nu\left( \sum_{(\phi,\mu)} \mu^\nu\frac{\ptl (f_a + f_bi
+  f_c j+  f_dk)}{\ptl q_\phi}
  \right) \notag \\
  & = \frac{1}{4}\sum_{\nu}dq^\nu  \left( \sum_{(\phi,\mu)} \mu^\nu\frac{\ptl f}{\ptl q_\phi}
\right),
\end{align}
which is the basis for the definitions for the right restricted HR derivatives
as given in the main text.

\section{Additional details for the Proof of Lemma \ref{lm:power} \label{app:lm1}}
To prove Lemma
\ref{lm:power}, we have used the following relation
\be
\frac{\ptl q^{-1}}{\ptl q} = - q^{-1}R(q^{-1}).
\ee
To show this result, we note $\ptl (q q^{-1})/\ptl q = \ptl 1/\ptl q = 0$.
Thus
\begin{align}
  0 & = q\frac{\ptl q^{-1}}{\ptl q} + \frac{1}{4}(q^{-1}-iq^{-1}i -
  jq^{-1}j-kq^{-1}k)             \notag\\
  &= q \frac{\ptl q^{-1}}{\ptl q} + R(q^{-1}),
\end{align}
from which the result follows. We have used equation (\ref{eq:d}) and the fact
that
\be
\frac{\ptl q}{\ptl q_a} = 1, \frac{\ptl q}{\ptl q_b} = i, \frac{\ptl q}{\ptl
q_c} = j, \frac{\ptl q}{\ptl q_d} = k.
\ee
The proof also uses the following recurrent relation
\be
\frac{\ptl q^{-n}}{\ptl q} = q^{-1}\left[\frac{\ptl q^{-(n-1)}}{\ptl q} -
R(q^{-n})\right],
\ee
which can be shown as follows: using the first product rule, we have
\begin{align}
  \frac{\ptl q^{-n}}{\ptl q} & = q^{-1}\frac{\ptl q^{-(n-1)}}{\ptl q} +
  \frac{1}{4} \left(\frac{\ptl q^{-1}}{\ptl q_a} q^{-{(n-1)}}
  -\frac{\ptl q^{-1}}{\ptl q_b} q^{-{(n-1)}}i\right.\notag\\
  &\phantom{=}
  \left.
    -\frac{\ptl q^{-1}}{\ptl q_c} q^{-{(n-1)}}j
    -\frac{\ptl q^{-1}}{\ptl q_d} q^{-{(n-1)}}k
  \right).
\end{align}
Using the fact $\ptl q q^{-1}/\ptl q_\phi = 0$ and the second product rule, we
can find
\be
\frac{\ptl q^{-1}}{\ptl q_\phi } = - q^{-1}\frac{\ptl q}{\ptl q_\phi} q^{-1}.
\ee
Thus
\begin{align}
  \frac{\ptl q^{-n}}{\ptl q} & = q^{-1}\frac{\ptl q^{-(n-1)}}{\ptl q}
  -\frac{q^{-1}}{4} \left(q^{-n}
  -iq^{-n}i\right.\notag\\
  &\phantom{=}
  \left.
    -jq^{-n}j
    -kq^{-n}k
  \right) \notag\\
  &=q^{-1}\frac{\ptl q^{-(n-1)}}{\ptl q} - q^{-1} R(q^{-n}).
\end{align}

\section{Derivations of the first chain rule \label{app:chain}}

   The function $f(g(q))$ may be view as a function of intermediate variables
   $g_a$, $g_b$, $g_c$ and $g_d$. Using the usual chain rule, we have
   \be
   \frac{\partial f}{\partial q_\beta} = \sum_{\phi} \frac{\ptl f}{\ptl
   g_\phi}\frac{\ptl g_\phi}{\ptl q_\beta},
   \ee
   with $\beta \in \{a,b,c,d\}$, which gives
   \be
   \nabla_r f = (\nabla_r^g f) P
   \ee
   where $P$ is a $4\times 4$ matrix with $P_{\phi \beta} = \ptl g_\phi / \ptl
   q_\beta$. With $(\nabla_r f) J^H = \nabla_q f$, and $\nabla_r^g f = 4
   (\nabla_q^g f) J$, the above equation leads to
   \be
   \nabla_q f = 4 (\nabla_q^g f) J P J^H,
   \ee
   where it is easy to show that $4 JP J^H = M$.


\end{document}